\documentclass[leqno]{amsart}
       
\usepackage{amsmath,amsfonts,amssymb}
\usepackage{amsthm}
\usepackage{hyperref}

\def\R{{\mathbb R}}
\def\N{{\mathbb N}}

\def\virgp{\raise 2pt\hbox{,}}

\def\bu{{\bf u}}

\def\({\left(}
\def\){\right)}
\def\le{\leqslant}
\def\ge{\geqslant}

\def\Eq#1#2{\mathop{\sim}\limits_{#1\rightarrow#2}}

\def\d{{\partial}}
\def\eps{\varepsilon}
\def\O{\mathcal O}

\theoremstyle{plain}
\newtheorem{theorem}{Theorem}[section]
\newtheorem{lemma}[theorem]{Lemma}
\newtheorem{corollary}[theorem]{Corollary}
\newtheorem{proposition}[theorem]{Proposition}

\theoremstyle{definition}

\theoremstyle{remark}
\newtheorem{remark}[theorem]{Remark}
\newtheorem*{remark*}{Remark}
\newtheorem*{example}{Example}

\numberwithin{equation}{section}

\begin{document}

\title[WKB for NLS and instability]{WKB analysis for the nonlinear
  Schr\"odinger equation and instability results}  
\author[R. Carles]{R{\'e}mi Carles}
\address{Wolfgang Pauli Institute, c/o Inst.~f.~Math.\\ Universit\"at Wien\\
        Nordbergstr.~15\\ A-1090 Wien\\ Austria\footnote{Future
  affiliation (starting September 2007): UMR CNRS 
  5149, Montpellier, France}}
\email{Remi.Carles@math.cnrs.fr}
\begin{abstract}
For the semi-classical limit of the cubic, defocusing
nonlinear Schr\"odinger equation with an external potential, we
explain the notion of criticality before a caustic is formed. In the
sub-critical and critical cases, we justify the WKB approximation.
In the super-critical case, the WKB analysis provides
a new phenomenon for the (classical) cubic, defocusing
nonlinear Schr\"odinger equation, which can be compared to the loss of
regularity established for the nonlinear wave equation by
G.~Lebeau. We also show some instabilities at the semi-classical level. 
\end{abstract}
\maketitle

\section{Introduction}
\label{sec:intro}

We consider the solution to the nonlinear Schr\"odinger equation (NLS)
with a cubic, defocusing nonlinearity, and an external potential:
\begin{align}
  i\eps \d_t u^\eps +\frac{\eps^2}{2}\Delta
  u^\eps &= V u^\eps + \eps^\kappa 
  |u^\eps|^2
  u^\eps,\quad (t,x)\in [0,T]\times \R^n,\label{eq:nlssemi}\\
u^\eps(0,x) &=
  a_0^\eps(x)e^{i\phi_0(x)/\eps} .\label{eq:CI}  
\end{align}
We are interested in the behavior of the solution $u^\eps$ as the
positive parameter $\eps$ goes to zero. According to the cultural
background, this field goes under the name of semi-classical limit,
WKB analysis (which is a little bit more specific, see below), or
geometrical optics. The general idea is to describe the asymptotic
behavior of $u^\eps$ with a simplified model, which involves geometric
quantities. In the case of \eqref{eq:nlssemi}, these quantities are
called either classical trajectories (in view of classical mechanics),
or rays (in view of geometric optics). 
\smallbreak

There are at least two motivations for such a
study. We outline them here, and refer to the survey \cite{RauchUtah}
for a broader discussion on this subject. The first one comes from the
applied mathematics, and may find 
its origins in physics. In the case of \eqref{eq:nlssemi}, suppose
that $\eps$ 
represents the (rescaled) Planck constant. It may be small compared to
the other parameters at stake. In this case, it is sensible to consider
that the asymptotic behavior of $u^\eps$ as $\eps\to 0$ provides a
reliable approximation of the exact solution. Hopefully, the
asymptotic model is easier to describe than the initial one
\eqref{eq:nlssemi}--\eqref{eq:CI}. Another motivation stems from the
propagation of singularities for equation where the small parameter
$\eps$ is not necessarily present initially. Most of the studies in
this direction concern hyperbolic equations. However, the belief
according to which Schr\"odinger equations share properties with
hyperbolic equations in the semi-classical limit is a first hint that
this field is applicable to Schr\"odinger equations as well (see
e.g. \cite{LebeauSchrod, JeremieAIF}). To illustrate this
statement, we give a result, whose proof will be straightforward after
the analysis of \eqref{eq:nlssemi}--\eqref{eq:CI}.
\begin{theorem}[\cite{CaARMA}, Cor.~1.7]\label{theo:flow}
   Let $n\ge 3$. Consider the cubic,
  defocusing NLS:
  \begin{equation}
    \label{eq:NLScubic}
    i\d_t u +\frac{1}{2}\Delta u = |u|^2 u\quad ; \quad u_{\mid
    t=0}=u_0\, .
  \end{equation}
Denote $s_c = \frac{n}{2}-1$. Let $0<s<s_c$. We can find a
  family $(u_0^\eps)_{0<\eps \le 
  1}$ in ${\mathcal S}({\mathbb R}^n)$ with 
\begin{equation*}
  \|u_0^\eps\|_{H^s({\mathbb R}^n)} \to 0 \text{ as }\eps \to 0\, ,
\end{equation*}
and $0<t^\eps \to 0$ such that the solution $u^\eps$ to
\eqref{eq:NLScubic} associated to $u_0^\eps$ satisfies: 
\begin{equation*}
  \|u^\eps(t^\eps)\|_{H^{k}({\mathbb R}^n)} \to +\infty \text{ as }\eps \to
 0\, , \ \forall k\in \left]\frac{s}{\frac{n}{2}-s}, s\right]\, .
\end{equation*}
\end{theorem}
This result is in the same spirit as the initial breakthrough by
G.~Lebeau (\cite{Lebeau01,Lebeau05}, see also
\cite{MetivierBourbaki}). These former 
results also rely on geometrical optics (in a super-critical
r\'egime; see below for this notion in the case of \eqref{eq:nlssemi}). For
\eqref{eq:NLScubic}, the above result was first established by
M.~Christ, J.~Colliander and T.~Tao \cite{CCT2} in the case $k=s$ (see
also \cite[Appendix]{BGTENS} and \cite[Appendix~B]{CaARMA}). The fact
that we can go strictly below the value $k=s$ stems from an analysis
of \eqref{eq:nlssemi} in a case where the nonlinearity should be
considered as quasilinear, and not semilinear. This result is then a
consequence of the original idea of E.~Grenier \cite{Grenier98}. 
\smallbreak

The rest of this text is organized as follows. In \S\ref{sec:formal},
we introduce the notion of criticality for \eqref{eq:nlssemi} at a
formal level. In \S\ref{sec:wkb}, we explain how to justify this
notion, and describe the asymptotic behavior of $u^\eps$ in different
cases. The proof of Theorem~\ref{theo:flow} is given in
\S\ref{sec:flow}. More instability results are given in
\S\ref{sec:instab}.

\section{WKB analysis and the notion of criticality}
\label{sec:formal}

This section remains at a formal level only. It is a preparation to
the forthcoming rigorous justifications. 
In a WKB analysis, one assumes for instance that the initial profile
$a_0^\eps$ can be expanded as a power series in $\eps$:
\begin{equation*}
  a_0^\eps \sim a_0 +\eps a_1 +\eps^2 a_2+\ldots
\end{equation*}
Note that the modulation factor $\eps^\kappa$ in front of the
nonlinearity in \eqref{eq:nlssemi} could be taken equal to one, up to
replacing $a_0^\eps$ with $\eps^{\kappa/2}a_0^\eps$. In particular, it
is clear that for large $\kappa>0$, the nonlinearity is expected to be
negligible in the limit $\eps\to 0$ (at least locally in time). WKB
analysis consists in seeking an approximation of the form
\begin{equation*}
  u^\eps(t,x) \Eq \eps 0 a(t,x)e^{\phi(t,x)/\eps},
\end{equation*}
where the amplitude $a$ and the phase $\phi$ are independent of
$\eps$. Plugging this approximate solution into \eqref{eq:nlssemi} and
canceling the first powers of $\eps$ yields: 
\begin{align}
  \O\(\eps^0\):\quad & \d_t \phi + \frac{1}{2}|\nabla \phi|^2 +V =0
  \ ;\ \phi_{\mid t=0}=\phi_0\quad\text{if }\kappa
  >0,\label{eq:eikonal}\\ 
\O\(\eps^{1}\):\quad & \d_t a +\nabla
  \phi\cdot \nabla a + 
  \frac{1}{2}a\Delta \phi =
\left\{
  \begin{aligned}
    &0 & \text{ if }\kappa >1,\\
    &-i|a|^2 a& \text{ if }\kappa =1.
  \end{aligned}
\right.\ ;\ a_{\mid t=0}=a_0\label{eq:transport}
\end{align}
This shows that the minimal value of $\kappa$ for which nonlinear
effects affect the solution at leading order is $\kappa=1$. We shall
consider here the following three cases: $\kappa >1$ (sub-critical
case), $\kappa =1$ (critical case), and $\kappa =0$ (a super-critical
case). We refer to \cite{CaBKW} for the case $0<\kappa<1$. The first
step consists in solving \eqref{eq:eikonal}. The potential $V$ may be
time dependent: $V=V(t,x)$.  
\begin{lemma}\label{lem:eikonal}
  Assume that $V$ and $\phi_0$ are smooth and sub-quadratic: 
  \begin{itemize}
  \item $V\in C^\infty(\R_t\times \R^n_x)$, and $\partial_x^\alpha V\in
  L^\infty_{\rm loc}(\R_t ;L^\infty(\R^n_x))$ as soon as $|\alpha|\ge 2$.
  \item $\phi_0\in C^\infty( \R^n)$, and $\partial_x^\alpha \phi_0\in
  L^\infty(\R^n)$ as soon as $|\alpha|\ge 2$.
  \end{itemize}
Then there exist $T>0$ and a unique
  solution $\phi_{\rm eik}\in C^\infty([0,T]\times\R^n)$ to
  \eqref{eq:eikonal}. 
This solution is sub-quadratic: $\partial_x^\alpha \phi_{\rm eik} \in
  L^\infty([0,T]\times\R^n)$ as soon as $|\alpha|\ge 2$. 
\end{lemma}
The proof of this lemma relies on Hamilton--Jacobi theory. Note that
the time of existence $T$ is uniform with respect to $x\in \R^n$: a
global inversion theorem is needed, which can be found in
\cite{SchwartzBook} or \cite{DG}. We refer to \cite{CaBKW} for the
proof of this lemma, and for a discussion on the optimality of the
assumptions. In particular, one should not expect the above solution
to remain smooth for all time. The appearance of singularities
corresponds to the formation of caustics. The aim of WKB analysis is
to describe the solution $u^\eps$ before a caustic is formed. See
e.g. \cite{CaIUMJ} for the description of a solution to
\eqref{eq:nlssemi}--\eqref{eq:CI} beyond a focal point. 
\smallbreak
To analyze \eqref{eq:transport}, we introduce the
Hamiltonian flow, on which the proof of Lemma~\ref{lem:eikonal}
relies: 
\begin{equation}
  \label{eq:hamilton}
\left\{
  \begin{aligned}
   &\partial_t x(t,y) = \xi \left(t,y\right)\quad ;\quad x(0,y)=y,\\ 
   &\partial_t \xi(t,y) = -\nabla_x V\left(t,x(t,y)\right)\quad ;\quad
   \xi(0,y)=\nabla \phi_0(y).
  \end{aligned}
\right.
\end{equation}
The time $T$ is such that the map $y\mapsto x(t,y)$ is a
diffeomorphism of $\R^n$ for $t\in [0,T]$. 
The key observation is that \eqref{eq:transport} is a
transport equation, which turns out to be an ordinary differential
equation along the classical trajectories. Introduce the Jacobi determinant
\begin{equation*}
  J_t(y) ={\rm det}\nabla_y x(t,y). 
\end{equation*}
Denote
\begin{equation*}
  A(t,y) := a \left(t, x(t,y)
  \right)\sqrt{J_t(y)}. 
\end{equation*}
For $t\in [0,T]$, \eqref{eq:transport}  is equivalent to:  
\begin{equation*}
  \d_t A = \left\{
\begin{aligned}
    &0 & \text{ if }\kappa >1,\\
    &-i J_t(y)^{-1}\left|A \right|^2
  A & \text{ if }\kappa =1.
  \end{aligned}
\right.
\quad ; \quad A(0,y)=a_0(y). 
\end{equation*}
This ordinary differential equation along the rays of geometrical
optics can be solved explicitly: we see that
$\partial_t |A|^2 =0$, hence, for $\kappa =1$,
\begin{equation*}
  A(t,y) = a_0(y) \exp\left(-i\int_0^t
  J_s(y)^{-1}\left|a_0(y) \right|^2 ds\right). 
\end{equation*}
Inverting the map $y\mapsto x(t,y)$ yields $a(t,x)$. We see that the
critical nonlinear effect is a self-modulation of the amplitude. In
the context of laser physics, this phenomenon is known as
\emph{phase self-modulation} (see e.g. \cite{ZS,Boyd,Donnat}). 
\smallbreak

We now turn to the super-critical case $\kappa=0$. To illustrate the
difficulty of this case, seek a more precise asymptotic expansion of
$u^\eps$: 
\begin{equation*}
  u^\eps(t,x)\sim \left({\tt a}_0(t,x) + \eps {\tt
  a}_1(t,x) + \eps^2 
  {\tt a}_2(t,x)+\ldots \right) e^{i\phi(t,x)/\eps}\, .
\end{equation*}
Plugging such an
asymptotic expansion into
\eqref{eq:nlssemi} yields a shifted cascade of equations:
\begin{align*}
 \O\left(\eps^0\right):&\quad \d_t \phi
 +\frac{1}{2}|\nabla\phi|^2 + V + |{\tt a}_0|^2=0\ ;\ \phi_{\mid
 t=0}=\phi_0.\\  
\O\left(\eps^1\right):&\quad \d_t {\tt a}_0 +\nabla\phi
 \cdot \nabla {\tt a}_0 +\frac{1}{2}{\tt a}_0\Delta \phi = 
-2i{\rm Re}\left({\tt a}_0\overline{{\tt
 a}_1}\right){\tt a}_0 \ ;\ {\tt a}_{0\mid
 t=0}=a_0.
\end{align*}
Two comments are in order. First, we see that there is a strong
coupling between the phase and the main amplitude: ${\tt a}_0$ is present
in the equation for $\phi$. Second, the above system is not closed:
$\phi$ is determined in function of ${\tt a}_0$, and ${\tt a}_0$ is
determined in function of ${\tt a}_1$. Even if we pursued the cascade
of equations, this phenomenon would remain: no matter how many terms
are computed, the system is never closed (see \cite{PGX93}). This is a
typical feature of super-critical cases in nonlinear geometrical optics
(see \cite{CheverryBullSMF,CG05}). 
\smallbreak

In the case when $V\equiv 0$ and $\phi_0\in H^s$, this problem was
resolved by E.~Grenier \cite{Grenier98}, by modifying the usual WKB
methods (see \S\ref{sec:super}). Note
that even though ${\tt a}_1$ is not determined by the above system, the
pair $(\rho,v):= (|{\tt a}_0|^2,\nabla \phi)$ solves a compressible Euler
equation:
\begin{equation}
  \label{eq:eulerV}
  \begin{aligned}
    \partial_t v +v\cdot \nabla v + \nabla V + \nabla \rho&=0\ ;\quad
    v\big|_{ t=0}=\nabla \phi_0\\
    \partial_t \rho+ \nabla\cdot (\rho v) &=0\ ;\quad \rho\big|_{t=0}=|a_0|^2.
  \end{aligned}
\end{equation}

\section{Rigorous WKB analysis}
\label{sec:wkb}
We outline here some results presented in details in \cite{CaBKW}. 

\subsection{Sub-critical and critical cases}
\label{sec:sous}

Let $\kappa \ge 1$. Under the assumptions of Lemma~\ref{lem:eikonal},
we change the Cauchy problem \eqref{eq:nlssemi}--\eqref{eq:CI}: define
$a^\eps$ by 
\begin{equation*}
  u^\eps (t,x)=a^\eps(t,x)e^{i\phi_{\rm eik}(t,x)/\eps}. 
\end{equation*}
Then for $t\in [0,T]$, \eqref{eq:nlssemi}--\eqref{eq:CI} is equivalent
to:
\begin{equation}
  \label{eq:asub}
  \begin{aligned}
    \d_t a^\eps +\nabla \phi_{\rm eik}\cdot \nabla
  a^\eps +\frac{1}{2}a^\eps\Delta \phi_{\rm eik}
  &=i\frac{\eps}{2}\Delta a^\eps 
  -i\eps^{\kappa -1}|a^\eps|^2 a^\eps ,\\
a^\eps_{\mid  t=0}&=a_0^\eps. 
  \end{aligned}
\end{equation}
Two things must be noticed: first, the potential $V$ and the initial phase
$\phi_0$ do not appear in this new problem. Second, the factors
involving $\phi_{\rm eik}$ have the following features: in view of
Lemma~\ref{lem:eikonal}, the term
$\Delta \phi_{\rm eik}$ is in $L^\infty([0,T]\times\R^n)$, and the
operator $\d_t+ \nabla \phi_{\rm eik}\cdot \nabla$ is a transport
operator. We can then obtain energy estimates in Sobolev spaces, and
establish:
\begin{proposition}\label{prop:asub}
  Under the assumptions of Lemma~\ref{lem:eikonal}, assume moreover that
  $a_0^\eps$ is bounded in $H^s(\R^n)$ uniformly for $\eps\in ]0,1]$,
  for all $s\ge 0$. Let
  $\kappa \ge 1$. Then for all $\eps\in ]0,1]$, \eqref{eq:asub} has a
  unique solution $a^\eps \in
  C^\infty([0,T]\times \R^n)\cap C([0,T];H^s)$ for all
  $s>n/2$. Moreover, $a^\eps$ is bounded in 
  $L^\infty([0,T]; H^s)$ uniformly in $\eps \in ]0,1]$, for all $s\ge
  0$. 
\end{proposition}
These uniform estimates allow us to neglect the term $\eps\Delta
a^\eps$ on the right hand side of \eqref{eq:asub}. Assume moreover
the following convergence: 
\begin{equation}\label{eq:cvci}
  a_0^\eps \to a_0 \text{ in }H^s(\R^n),\quad \forall s\ge 0. 
\end{equation}
\begin{corollary}\label{cor:simplifsub}
 Let $\kappa \ge 1$. Under the above assumptions, 
  \begin{equation*}
  \left\| a^\eps -   \widetilde a^\eps
  \right\|_{L^\infty([0,T];H^s)}\to 0 \quad \text{as }\eps \to
  0,\quad \forall s\ge 0, 
  \end{equation*}
where $\widetilde a^\eps$ solves:
\begin{equation}\label{eq:atildeeps}
  \d_t \widetilde a^\eps +\nabla \phi_{\rm eik}\cdot
  \nabla \widetilde a^\eps + 
  \frac{1}{2}\widetilde a^\eps\Delta \phi_{\rm eik} =
  -i\eps^{\kappa -1}|\widetilde
  a^\eps|^2\widetilde a^\eps \quad ; \quad 
\widetilde a^\eps_{\mid  t=0}=a_0.
\end{equation}
\end{corollary}
Proceeding as in \S\ref{sec:formal}, 
denote
\begin{equation*}
  A^\eps(t,y) := \widetilde a^\eps \left(t, x(t,y)
  \right)\sqrt{J_t(y)}. 
\end{equation*}
We see that so long as $y\mapsto x(t,y)$ defines a global
diffeomorphism (which is guaranteed for $t\in [0,T]$ by construction),
\eqref{eq:atildeeps} is equivalent to:  
\begin{equation*}
  \partial_t A^\eps = -i\eps^{\kappa
  -1}J_t(y)^{-1}\left|A^\eps \right|^2
  A^\eps\quad ; \quad A^\eps(0,y)=a_0(y). 
\end{equation*}
This ordinary differential equation along the rays of geometrical
optics can be solved explicitly: 
\begin{equation*}
  A^\eps(t,y) = a_0(y) \exp\left(-i\eps^{\kappa -1} \int_0^t
  J_s(y)^{-1}\left|a_0(y) \right|^2 ds\right). 
\end{equation*}
Back to the initial solution $u^\eps$, we conclude:
\begin{proposition}\label{prop:sub}
  Let
  $\kappa \ge 1$. Under the above assumptions, for all $\eps \in ]0,1]$,
  \eqref{eq:nlssemi}-\eqref{eq:CI} has a unique solution $u^\eps \in
  C^\infty([0,T]\times \R^n)\cap C([0,T];H^s)$ for all $s>n/2$. Moreover,
  there exist $ a, G\in C^\infty([0,T]\times \R^n)$,
  independent of $\eps \in ]0,1]$, where
  $ a \in C([0,T]; L^2\cap L^\infty)$, and $G$ is
  real-valued with $G \in C([0,T];   L^\infty)$, such that:
  \begin{equation*}
    \left\| u^\eps -  a e^{i\eps^{\kappa -1}G}e^{i\phi_{\rm eik}
        /\eps}\right\|_{L^\infty([0,T]; L^2\cap L^\infty) } \to 0\quad
    \text{as }\eps \to 0.
  \end{equation*}
The profile $ a$ solves the initial value problem:
\begin{equation}\label{eq:alibre}
  \partial_t  a +\nabla \phi_{\rm eik}\cdot \nabla a +
  \frac{1}{2} a\Delta \phi_{\rm eik} =0\quad ;\quad a_{\mid
    t=0}=a_0, 
\end{equation}
and $G$  depends nonlinearly on $ a$:
\begin{equation*}
  \begin{aligned}
    a(t,x) &= \frac{1}{\sqrt{J_t(y(t,x))}}a_0\left(y(t,x)\right),\\
  G(t,x) &= -\int_0^t
  J_s(y(t,x))^{-1}\left|a_0(y(t,x)) \right|^2 ds.
  \end{aligned}
\end{equation*}
In particular, if $\kappa>1$, then 
\begin{equation*}
    \left\| u^\eps - a e^{i\phi_{\rm eik}
        /\eps}\right\|_{L^\infty([0,T]; L^2\cap L^\infty) } \to 0\quad
    \text{as }\eps \to 0,
  \end{equation*}
and no nonlinear effect is present in the leading order behavior of
$u^\eps$. If $\kappa =1$, nonlinear effects are present at leading
order, measured by $G$.
\end{proposition}

\subsection{Super-critical case: $\kappa =0$}
\label{sec:super}
In this case, we recall the beautiful idea of E.~Grenier, which makes
it possible to consider the case $V\equiv 0$ and $\phi_0\in H^s(\R^n)$
for sufficiently large $s$. The approach consists in reversing the
steps of the WKB analysis: usually, one seeks an approximate solution,
and then tries to show that stability arguments imply that the exact
solution is well approximated by this process. To overcome the issues
mentioned in \S\ref{sec:formal}, the idea in \cite{Grenier98}
consists in first writing the unknown as:
\begin{equation*}
  u^\eps(t,x)=a^\eps(t,x)e^{i\phi^\eps(t,x)/\eps},
\end{equation*}
where the amplitude $a^\eps$ is \emph{complex-valued} (even if
$a_0^\eps$ is real-valued), and $\phi^\eps$ is real-valued. Doing so,
one introduces one degree of freedom to rewrite
\eqref{eq:nlssemi}--\eqref{eq:CI}. The usual approach in physics (see
e.g. \cite{LandauQ})
consists in writing
\begin{equation*}
\left\{
  \begin{aligned}
   \d_t \phi^\eps +\frac{1}{2}\left|\nabla \phi^\eps\right|^2 + 
    |a^\eps|^2= \eps^2 \frac{\Delta a^\eps}{2
    a^\eps}\quad &; \quad 
    \phi^\eps\big|_{t=0}=\phi_0\, ,\\
\d_t a^\eps +\nabla \phi^\eps \cdot \nabla a^\eps +\frac{1}{2}a^\eps
\Delta \phi^\eps  = 0\quad  &;\quad
a^\eps\big|_{t=0}= a^\eps_0 . 
  \end{aligned}
\right.     
\end{equation*}
Obviously, this approach is delicate when $a^\eps$ has zeroes; see the
discussion in \cite{PGX93} on this subject. The choice in
\cite{Grenier98} is to write:
\begin{equation}\label{eq:systexact0}
\left\{
  \begin{aligned}
    \d_t \phi^\eps +\frac{1}{2}\left|\nabla
    \phi^\eps\right|^2 + 
    |a^\eps|^2= 0\quad &; \quad
    \phi^\eps\big|_{t=0}=\phi_0\, ,\\
\d_t a^\eps +\nabla \phi^\eps \cdot \nabla
    a^\eps +\frac{1}{2}a^\eps 
\Delta \phi^\eps  = i\frac{\eps}{2}\Delta
    a^\eps\quad & ;\quad a^\eps\big|_{t=0}= a^\eps_0\, . 
  \end{aligned}
\right.
\end{equation}
Inspired by the fact that the expected limit is related to the
compressible Euler equation (see \S\ref{sec:formal}), introduce the
``velocity'' $v^\eps = \nabla \phi^\eps$. Then \eqref{eq:systexact0}
yields: 
\begin{equation}\label{eq:systexact}
\left\{
  \begin{aligned}
    \d_t v^\eps +v^\eps \cdot \nabla v^\eps
    + 2 
    \operatorname{Re}\left(\overline{a^\eps}\nabla
    a^\eps\right)= 
    0\quad &; \quad 
    v^\eps\big|_{t=0}=\nabla \phi_0\, ,\\
\d_t a^\eps +v^\eps\cdot \nabla a^\eps
    +\frac{1}{2}a^\eps 
\nabla\cdot v^\eps  = i\frac{\eps}{2}\Delta
    a^\eps\quad & ;\quad 
a^\eps\big|_{t=0}= a^\eps_0\, . 
  \end{aligned}
\right.
\end{equation}
Separate real and imaginary parts of $a^\eps$, $a^\eps = 
a_1^\eps + ia_2^\eps$. Then we have 
\begin{equation}
  \label{eq:systhyp}
  \partial_t \bu^\eps +\sum_{j=1}^n
  A_j(\bu^\eps)\partial_j \bu^\eps 
  = \frac{\eps}{2} L 
  \bu^\eps\, , 
\end{equation}
\begin{equation*}
  \text{with}\quad \bu^\eps = \left(
    \begin{array}[l]{c}
       a_1^\eps \\
       a_2^\eps \\
       v^\eps_1 \\
      \vdots \\
       v^\eps_n
    \end{array}
\right)\quad , \quad L = \left(
  \begin{array}[l]{ccccc}
   0  &-\Delta &0& \dots & 0   \\
   \Delta  & 0 &0& \dots & 0  \\
   0& 0 &&0_{n\times n}& \\
   \end{array}
\right),
\end{equation*}
\begin{equation*}
  \text{and}\quad A(\bu,\xi)=\sum_{j=1}^n A_j(\bu)\xi_j
= \left(
    \begin{array}[l]{ccc}
      v\cdot \xi & 0& \frac{a_1 }{2}\,^{t}\xi \\ 
     0 &  v\cdot \xi & \frac{a_2}{2}\,^{t}\xi \\ 
     2  a_1 \, \xi
     &2  a_2\, \xi &  v\cdot \xi I_n 
    \end{array}
\right).
\end{equation*}
The matrix $A(\bu,\xi)$ can be symmetrized by 
\begin{equation*}
  S=\left(
    \begin{array}[l]{cc}
     I_2 & 0\\
     0& \frac{1}{4}I_n
    \end{array}
\right).
\end{equation*}
The important point to notice is that the operator $L$ is
skew-symmetric: it is invisible in the energy estimates, so the loss
of derivative (it is the only operator of order two in
\eqref{eq:systhyp}) is avoided. Denote $H^\infty =\cap_{s\ge
  0}H^s(\R^n)$. 
Classical theory on symmetric hyperbolic systems yields a solution
$(v^\eps,a^\eps)$ to \eqref{eq:systexact}.
Once $v^\eps$ is known, we note that it is irrotational, so
there exists $\phi^\eps$ such that $v^\eps=\nabla \phi^\eps$. Up to
adding a function of time only, $(\phi^\eps,a^\eps)$ solves
\eqref{eq:systexact0}. 
\begin{proposition}[\cite{Grenier98}, Th.~1.1]\label{prop:p}
Let $\kappa =0$. Suppose that $\phi_0\in H^\infty$, and that
  $a_0^\eps$ is bounded in 
  $H^s(\R^n)$ uniformly for $\eps\in ]0,1]$, for all $s\ge 0$. Let
$s>2+n/2$. There exist $T_s>0$ independent of $\eps\in ]0,1]$  and
$u^\eps = a^\eps e^{i\phi^\eps/\eps}$
solution to \eqref{eq:nlssemi}--\eqref{eq:CI}   on
  $[0,T_s]$. Moreover, $a^\eps$ and 
$\phi^\eps$ are bounded in $L^\infty([0,T_s];H^s)$,
uniformly in $\eps\in ]0,1]$. 
\end{proposition}
Assume moreover that \eqref{eq:cvci} holds. 
The solution to \eqref{eq:systexact0} formally converges to the
solution  of:
\begin{equation}\label{eq:systlim}
\left\{
  \begin{aligned}
    \d_t \phi +\frac{1}{2}\left|\nabla \phi\right|^2 + 
    |a|^2= 0\quad &; \quad
    \phi\big|_{t=0}=\phi_0\, ,\\
\d_t a +\nabla \phi \cdot \nabla a +\frac{1}{2}a
\Delta \phi  = 0\quad & ;\quad
a\big|_{t=0}= a_0\, . 
  \end{aligned}
\right.
\end{equation}
Under the above assumptions, 
\eqref{eq:systlim} has a unique solution $(a,\phi)\in
L^\infty([0,T_*];H^m)^2$ for all $m>0$ for some $T_*>0$ 
independent of $m$ (see e.g. \cite{AlinhacGerard,Majda}). 
\begin{remark}\label{rem:nl}
  More general nonlinearity. Suppose that we consider a more general
  nonlinearity:
  \begin{equation*}
    i\eps \d_t u^\eps +\frac{\eps^2}{2}\Delta u^\eps = f\(
    |u^\eps|^2\)u^\eps. 
  \end{equation*}
Then following the same lines as above, the symmetrizer naturally
becomes
\begin{equation*}
  S=\left(
    \begin{array}[l]{cc}
     I_2 & 0\\
     0& \frac{1}{4f'(|a^\eps|^2)}I_n
    \end{array}
\right).
\end{equation*}
For this matrix to be positive, and to be able to estimate its time
derivative, it is natural to assume $f'>0$. This corresponds to the
assumption made in \cite{Grenier98}, and in \cite{CaBKW}. For the
above analysis to be 
valid, the nonlinearity has to be defocusing, and cubic at the
origin. In particular, the WKB analysis for the quintic defocusing NLS
is still an open problem. Note however that it is possible to
construct solutions to the limit problem in that case (the analogue of
\eqref{eq:systlim}), thanks to the result of \cite{MUK86} and the
geometrical analysis of \S\ref{sec:sous}. Yet, the nonlinear change
of variable of \cite{MUK86} is apparently incompatible with the above
remark that $L$ is skew-symmetric.
\end{remark}

If we
suppose in addition that
there exists $a_0,a_1 \in H^\infty$
  such that
  \begin{equation}\label{eq:10h27}
    a_0^\eps = a_0 +\eps a_1 +o(\eps)\quad
    \text{in }H^s, \ \forall s\ge 0,
  \end{equation}
then we infer: 
\begin{proposition}\label{prop:estprec}
Let $s\in \N$. Then $T_s\ge T_*$, and there exists $C_s$ independent of
$\eps$ such that for 
every $0\le t\le T_*$, 
\begin{equation*}
  \| a^\eps (t)- a(t)\|_{H^s}\le C_s \eps\quad ;\quad
  \| \phi^\eps(t) - \phi(t)\|_{H^s} \le C_s \eps t . 
\end{equation*}
\end{proposition}
Note that this suffices to describe $u^\eps$ for very small time only:
\begin{align*}
  u^\eps - ae^{i\phi/\eps}& = a^\eps e^{i\phi^\eps /\eps} -
  ae^{i\phi/\eps} \\
&= \(a^\eps -a \) e^{i\phi^\eps /\eps} +2i ae^{i(\phi^\eps +\phi^\eps)/2\eps} 
 \sin \(\frac{\phi^\eps -\phi}{2\eps}\).
\end{align*}
The first term of the right hand side is of order $\O(\eps)$ in
$L^2\cap L^\infty$, but the second one is of order $\O(t)$ only:
therefore, we only have
\begin{equation*}
  u^\eps(t,x)\sim a(t,x)e^{i\phi(t,x)/\eps}\quad \text{ for }0\le
  t\ll 1. 
\end{equation*}
To have a better error estimate, it is necessary to compute the next
term in the asymptotic expansion of $(\phi^\eps,a^\eps)$ in powers of
$\eps$. For times of order $\O(1)$, the initial corrector $a_1$ must
be taken into account: 
\begin{proposition}\label{prop:correc}
  Define
  $(a^{(1)},\phi^{(1)})$ by 
\begin{equation*}
  \begin{aligned}
    \d_t \phi^{(1)} +\nabla \phi \cdot \nabla \phi^{(1)} +
    2\operatorname{Re}\left(\overline a a^{(1)}\right)&=0,\\ 
   \d_t a^{(1)} +\nabla\phi\cdot \nabla a^{(1)} + \nabla
   \phi^{(1)}\cdot \nabla a + \frac{1}{2} a^{(1)}\Delta \phi
   +\frac{1}{2}a\Delta \phi^{(1)}      &= \frac{i}{2}\Delta a,\\
\phi^{(1)}\big|_{t=0}=0\quad ; \quad a^{(1)}\big|_{t=0}=a_1.
  \end{aligned}
\end{equation*}
Then $a^{(1)},\phi^{(1)}\in
L^\infty([0,T_*];H^s)$ for every $s\ge 0$, and
\begin{equation*}
  \|a^\eps - a - \eps a^{(1)}\|_{L^\infty([0,T_*];H^s)}+
  \|\phi^\eps - \phi - \eps 
  \phi^{(1)}\|_{L^\infty([0,T_*];H^s)} \le C_s\(\eps^2+o(\eps)\),\quad
  \forall s\ge 0\, . 
\end{equation*}
\end{proposition}
Despite 
the notations, it seems unadapted to consider $\phi^{(1)}$ as being
part of the phase. Indeed, we infer from Proposition~\ref{prop:correc}
that 
\begin{equation*}
  \left\|u^\eps - a e^{i\phi^{(1)}}
    e^{i\phi/\eps}\right\|_{L^\infty([0,T_*];L^2\cap
    L^\infty)}= o(1).  
\end{equation*}
\begin{remark}
 If the term $o(\eps)$ in \eqref{eq:10h27} is controlled more
precisely as a $\O(\eps^2)$, then the above $o(1)$ becomes a
$\O(\eps)$.  
\end{remark}
Relating this information to the WKB methods presented at the end of
\S\ref{sec:formal}, we would have: 
\begin{equation*}
  {\tt a}_0 = a e^{i\phi^{(1)}}.
\end{equation*}
Since $\phi^{(1)}$ depends on $a_1$ while $a$ does not, we retrieve
the fact that in super-critical r\'egimes, the leading order amplitude
in WKB methods depends on the initial first corrector $a_1$. 
\begin{remark}
  The term $e^{i\phi^{(1)}}$ does not appear in the Wigner measure
  of $a e^{i\phi^{(1)}} e^{i\phi/\eps}$. Thus, from the point of view of
    Wigner measures, the asymptotic behavior of the exact solution is
    described by the Euler-type system \eqref{eq:eulerV}. 
\end{remark}
\begin{remark}
  If we assume that $a_0$ is real-valued, then so is $a$. If moreover
  $a_1$ is purely imaginary (for instance, if $a_1=0$), then we see
  that $a^{(1)}$ is purely imaginary, hence, $\phi^{(1)}\equiv
  0$. 
\end{remark}
So far we have assumed $V\equiv 0$ and $\phi_0\in H^\infty$. If we try
to mimic the approach of \cite{Grenier98} for a non-trivial external
potential for instance, we have to consider:
\begin{align*}
    \partial_t \phi^\eps +\frac{1}{2}\left|\nabla
    \phi^\eps\right|^2 + V+ 
    |a^\eps|^2= 0\quad &; \quad
    \phi^\eps\big|_{t=0}=\phi_0\, ,\\
\partial_t a^\eps +\nabla \phi^\eps \cdot \nabla
    a^\eps +\frac{1}{2}a^\eps 
\Delta \phi^\eps  = i\frac{\eps}{2}\Delta
    a^\eps\quad & ;\quad 
a^\eps\big|_{t=0}= a^\eps_0\, . 
\end{align*}
The analysis of \cite{Grenier98} works in the same way
only when $\nabla V\in L^\infty_{\rm loc}(\R_t;H^s(\R^n))$ for a
sufficiently large $s$. To be able to consider general sub-quadratic
potentials (including the harmonic oscillator), resume the assumption
of Lemma~\ref{lem:eikonal}, and write
\begin{equation*}
  \phi^\eps =\phi_{\rm eik}+\varphi^\eps.
\end{equation*}
Working with the unknown $(\varphi^\eps,a^\eps)$, we see that we are
now rid of the external potential $V$, and of the possibly unbounded
initial phase $\phi_0$. The price to pay is that extra terms have
appeared. The good news however is that these extra terms are
\emph{semilinear} (as in \S\ref{sec:sous}), and can be treated by
perturbative methods in energy estimates.  We conclude:
\begin{theorem}\label{theo:BKWV}
  Let $\kappa =0$. Under the above assumptions, there exists $T_*>0$
  independent 
  of $\eps \in ]0,1]$ and a unique solution $u^\eps \in
  C^\infty([0,T_*]\times \R^n)\cap C([0,T_*];H^s)$ for all $s>n/2$ to
  \eqref{eq:nlssemi}--\eqref{eq:CI}. Moreover, there exist $a,\varphi \in
  C([0,T_*];H^s)$ for every $s\ge 0$, such that:
  \begin{equation*}
    \limsup_{\eps \to 0}\left\| u^\eps - a
    e^{i(\varphi+\phi_{\rm eik})/\eps}\right\|_{L^2\cap
    L^\infty}=\O(t)\quad \text{as }t\to 0. 
  \end{equation*}
Here, $a$ and $\varphi$ are nonlinear functions of $\phi_{\rm eik}$ and
$a_0$. Finally, there exists $\varphi^{(1)}\in
  C([0,T_*];H^s)$ for every $s\ge 0$, real-valued, such that: 
  \begin{equation*}
    \limsup_{\eps \to 0}\sup_{0\le t\le T_*}\left\|
    u^\eps - ae^{i\varphi^{(1)}} 
    e^{i(\varphi+\phi_{\rm eik})/\eps}\right\|_{L^2\cap L^\infty}=0. 
  \end{equation*}
The phase shift $\varphi^{(1)}$ is a nonlinear function of $\phi_{\rm
  eik}, a_0$ and $a_1$.  
\end{theorem}
\begin{remark}
  In \cite{CaBKW}, some assumptions on the momentum of $a_0^\eps$ are
  made, and not here. This is due to the fact that here, the
  nonlinearity that we consider is exactly cubic. When it is cubic at
  the origin only (see Remark~\ref{rem:nl}), extra estimates are
  needed, which apparently impose some extra decay at infinity for
  $a_0^\eps, a_0$ and $a_1$. 
\end{remark}

\section{Proof of Theorem~\ref{theo:flow}}
\label{sec:flow}

Theorem~\ref{theo:flow} is a straightforward consequence of
Proposition~\ref{prop:correc}.  For $a_0\in {\mathcal S}({\mathbb
  R}^n)$, let 
  \begin{equation*}
    u_0(x) = \lambda^{-\frac{n}{2}+s}a_0\left(\frac{x}{\lambda}\right) .
  \end{equation*}
Let $\eps = \lambda^{\frac{n}{2}-1 -s}$: $\eps$ and $\lambda$ go
simultaneously to 
zero, since $s<s_c$. Define 
\begin{equation*}
  \psi^\eps(t,x) = u^\lambda ( \eps t,x) =\lambda^{\frac{n}{2}-s}u\left(
  \lambda^{\frac{n}{2}+1-s} t,\lambda x\right)\, . 
\end{equation*}
It solves:
\begin{equation}\label{eq:psi52}
  i\eps\d_t \psi^\eps +\frac{\eps^2}{2}\Delta \psi^\eps =
  |\psi^\eps|^{2}\psi^\eps \quad 
  ;\quad \psi^\eps_{\mid t=0} = a_0(x) \, .
\end{equation}
The idea of the proof is that for times of order ${\mathcal O}(1)$,
$\psi^\eps$ has become $\eps$-oscillatory. 

We infer from Proposition~\ref{prop:correc} that there exist $T>0$
independent of $\eps\in ]0,1]$, and $a,\phi,\phi_1 \in C([0,T];H^m)$ for
any $m\ge 0$, such that:
\begin{equation*}
  \left\| \psi^\eps - a e^{i\phi_1} e^{i\phi/ \eps}\right\|_{L^\infty([0,T];
  H^m)} \le C_m \eps^{1-m}.
\end{equation*}
Since the $\dot H^m$-norm of $a e^{i\phi_1} e^{i\phi/ \eps}$ is of order
$\eps^{-m}$ (when $\phi$ is not stationary), we deduce that there exists
$t\in ]0,T]$ such that for any 
$m\ge 0$:
\begin{equation*}
  \| \psi^\eps(t)\|_{\dot H^m} \approx \eps^{-m}.
\end{equation*}
This implies: 
\begin{equation*}
  \left\| u \left( \lambda^{\frac{n}{2}+1 -s}t\right)\right\|_{\dot
  H^k} \approx 
  \lambda^{s-k}\| \psi^\eps(t)\|_{\dot H^k}\approx \lambda^{s-k} \eps^{-k} =
  \lambda^{s-k-k\left(\frac{n}{2}-1 
  -s\right)} \, . 
\end{equation*}
The result then follows when considering the limit $\lambda \to 0$. We
get exactly the statement of the theorem by 
replacing $a_0$ by $|\log \lambda|^{-1} a_0$ for instance. 
\begin{remark}\label{rem:laplace}
  The proof of ill-posed presented in \cite{CCT2} (see also
  \cite[Appendix]{BGTENS}, \cite[Appendix~B]{CaARMA}) consists in
  neglecting the Laplacian in \eqref{eq:psi52} for very small times,
  and integrating explicitly an ordinary differential equation. Proving
  that the Laplacian is negligible stems from Gronwall
  lemma. Essentially, the error satisfies an inequality of the form
  \begin{equation*}
    \|w^\eps(t)\|_X \lesssim \eps +\frac{1}{\eps}\int_0^t \|w^\eps(s)\|_Xds, 
  \end{equation*}
for some space $X$ that we do not describe. The singular factor
$\eps^{-1}$ is due to the $\eps$ in front of the time derivative, and
to the fact that no power of $\eps$ is present in front of the
nonlinearity. Therefore, Gronwall lemma yields no better than:
\begin{equation*}
  \|w^\eps(t)\|_X \lesssim \eps e^{Ct/\eps},
\end{equation*}
for some $C>0$. The error is small on an interval of the form $[0,\eps
|\log \eps|^\theta]$ for some $\theta>0$. This is enough to prove
Theorem~\ref{theo:flow} for $k=s$. This analysis considers
\eqref{eq:psi52} as a semilinear equation, since the nonlinearity is
viewed as a perturbation of the linear equation. To prove
Theorem~\ref{theo:flow} for $k<s$, it seems necessary to consider
\eqref{eq:psi52} as a quasilinear equation, as was done by E.~Grenier.
Note also that the quasilinear approach shows that the Laplacian in
\eqref{eq:psi52}  is negligible for $0<t^\eps \ll \eps^{1/3}$, that is
a ``much larger'' interval than $[0,\eps
|\log \eps|^\theta]$ (but still very small!). See the next section. 
\end{remark}
\begin{remark}
  On the other hand, Theorem~\ref{theo:flow} is valid only for cubic,
  defocusing nonlinear Schr\"odinger equations, while the results in
  \cite{CCT2} are valid for more general equations. This is due to the
  fact that the justification of super-critical nonlinear geometric
  optics for times $\O(1)$ is available only for nonlinearities which
  are defocusing, 
  and cubic at the origin (see Remark~\ref{rem:nl}). Since the proofs
  of ill-posedness rely on an homogeneous change of unknown function,
  we are left with the only possibility of an exactly cubic,
  defocusing nonlinearity. However, it is very likely that (an analogue
  of) Theorem~\ref{theo:flow} should be true under more general
  assumptions. 
\end{remark}
\section{Instability for the semi-classical equation}
\label{sec:instab}

The results we present in the paragraph are taken from
\cite{CaARMA}. We assume $V=\phi_0=0$ for the sake of concision. 
We first fix some notations. \\
{\bf Notation.}
Let $(\alpha^\eps)_{0<\eps\le 1}$ and $(\beta^\eps)_{0<\eps\le 1}$ be
two families 
of positive real numbers. 
\begin{itemize}
\item We write $\alpha^\eps \ll \beta^\eps$ if
$\displaystyle \limsup_{\eps\to 0}\alpha^\eps/\beta^\eps =0$.
\item We write $\alpha^\eps \lesssim \beta^\eps$ if 
$\displaystyle \limsup_{\eps\to 0}\alpha^\eps/\beta^\eps <\infty$.
\item We write $\alpha^\eps \approx \beta^\eps$ if $\alpha^\eps \lesssim
  \beta^\eps$ and $\beta^\eps \lesssim \alpha^\eps$. 
\end{itemize}
A typical result of \cite{CaARMA} is the following:
\begin{theorem}\label{theo:strong}
Let $n\ge 1$, $a_0 ,\widetilde a_0^\eps\in{\mathcal S}({\mathbb
  R}^n)$, where $a_0$ is independent of $\eps$.  Let $u^\eps$ and  
$v^\eps$ solve the initial value problems:
\begin{align*}
i\eps \d_t u^\eps + \frac{\eps^2}{2}\Delta u^\eps &=
|u^\eps|^2 u^\eps 
\ ; \ u^\eps\big|_{t=0}= a_0\, .\\
i\eps \d_t v^\eps + \frac{\eps^2}{2}\Delta v^\eps &= 
|v^\eps|^2 v^\eps 
\ ; \ v^\eps\big|_{t=0}= \widetilde a_0^\eps\, .
\end{align*}
Assume that there exists $N\in {\mathbb N}$ and $\eps^{1-\frac{1}{N}}\ll
\delta^\eps \ll 1$ such that:
\begin{equation}\label{eq:hyppert}
\begin{aligned}
  \left\|a_0 -\widetilde
    a_0^\eps\right\|_{H^s} \approx \delta^\eps\, ,\ 
  \forall s\ge 0\  ;\ 
\limsup_{\eps\to 0}\left\| \frac{\operatorname{Re}(a_0-\widetilde
    a_0^\eps)\overline{a_0}}{\delta^\eps}  \right\|_{L^\infty({\mathbb
    R}^n)} \not =0. 
\end{aligned}
\end{equation}
Then we can find $0<t^\eps\ll 1$ such that:
$\displaystyle \left\| u^\eps(t^\eps) - v^\eps(t^\eps) \right\|_{L^2}\gtrsim
1$. More precisely, this mechanism occurs as soon as $t^\eps \delta^\eps
\gtrsim \eps$. 
In particular, for all $s\ge 0$, 
\begin{equation*}
\frac{\left\| u^\eps - v^\eps
  \right\|_{L^\infty([0,t^\eps];L^2)}}{\left\| u^\eps_{\mid t=0} - v^\eps_{\mid
  t=0}
  \right\|_{H^s}}\to +\infty \quad \text{as }\eps\to 0\, .
\end{equation*}
\end{theorem}
\begin{example}
  Consider $a_0,b_0\in {\mathcal S}({\mathbb R}^n)$ independent of
  $\eps$, such that 
  $\operatorname{Re}(\overline{a_0}b_0)\not \equiv 0$, and take $\widetilde
  a_0^\eps = a_0  +\delta^\eps b_0$. 
\end{example}
\begin{example}
  Consider $a_0\in {\mathcal S}({\mathbb R}^n)$ independent of $\eps$ and
  $x^\eps\in {\mathbb R}^n$. We can take $\widetilde  a_0^\eps(x) = a_0
  (x-x^\eps)$, provided that $|x^\eps|=\delta^\eps$ and
\begin{equation*}
  \limsup_{\eps\to 0} \left\| \frac{x^\eps}{|x^\eps|}\cdot \nabla\left(
    |a_0|^2\right)\right\|_{L^\infty} \not =0. 
\end{equation*}
This example and the analysis of \cite{CaARMA} make it possible to
refine some results of \cite{BZ}. 
\end{example}
The general idea consists in using the WKB analysis in this
super-critical case. Roughly speaking, we have seen that
\eqref{eq:systlim} provides a good approximation of $u^\eps$ for very small
time only. Since we are interested in instabilities occurring for very
small time, this is not a problem for us now. The coupling in
\eqref{eq:systlim} shows that a small perturbation of $a_0$ yields a
small perturbation of $\phi$. But when we write
\begin{equation*}
  u^\eps(t,x)\sim a(t,x)e^{i\phi(t,x)/\eps},
\end{equation*}
we see that this small perturbation is divided by $\eps$, which goes
to zero. The result may not be small\ldots
\smallbreak

Technically, our approach consists in resuming the result provided by
Proposition~\ref{prop:p}. Instead of letting $\eps \to 0$ in the
initial data of \eqref{eq:systexact0}, just neglect the skew-symmetric
term (recall that we assume $\phi_0=0$):
\begin{equation}\label{eq:systexact1}
  \begin{aligned}
    \d_t \Phi^\eps +\frac{1}{2}\left|\nabla
    \Phi^\eps\right|^2 + 
    |{\bf a}^\eps|^2= 0\quad &; \quad
    \Phi^\eps\big|_{t=0}=0\, ,\\
\d_t {\bf a}^\eps +\nabla \Phi^\eps \cdot \nabla
    {\bf a}^\eps +\frac{1}{2}{\bf a}^\eps 
\Delta \Phi^\eps  = 0\quad & ;\quad {\bf a}^\eps\big|_{t=0}= a^\eps_0\, . 
  \end{aligned}
\end{equation}
Assuming that $a_0^\eps$ is bounded in $H^s$ for all $s\ge 0$, 
\eqref{eq:systexact1} has a unique solution $(\Phi^\eps,{\bf a}^\eps)\in
L^\infty([0,T_*];H^m)^2$ for all $m>0$ for some $T_*>0$ 
independent of $\eps$ and $m$.
\begin{proposition}\label{prop:estprec1}
Let $s\in \N$. Then $T_s\ge T_*$, and there exists $C_s$ independent of
$\eps$ such that for 
every $0\le t\le T_*$, 
\begin{equation*}
  \| a^\eps (t)- {\bf a}^\eps(t)\|_{H^s}\le C_s \eps t\quad ;\quad
  \| \phi^\eps(t) - \Phi^\eps(t)\|_{H^s} \le C_s \eps t^2 . 
\end{equation*}
\end{proposition}
The second idea consists in considering the Taylor expansion in time
of $(\Phi^\eps,{\bf a}^\eps)$:
\begin{equation*}
 \Phi^\eps(t,x)\sim \sum_{j\ge 1}t^{2j-1}\Phi_j^\eps(x)\quad ;\quad
 {\bf a}^\eps(t,x)\sim \sum_{j\ge 1}t^{2j}{\bf a}_j^\eps(x).  
\end{equation*}
Note that only odd powers of $t$ are present in the expansion of
$\Phi^\eps$, and even powers in that of ${\bf a}^\eps$. This is
because we have assumed $\phi_0=0$. Plugging these expansions into
\eqref{eq:systexact1}, we get formally:
\begin{equation*}
 {\bf a}_0^\eps=a_0^\eps\quad ;\quad \Phi_1^\eps = -|a_0^\eps|^2.
\end{equation*}
We can then check that a perturbation of order $\delta^\eps$ of
$a_0^\eps$ yields a perturbation of order $\delta^\eps$ of
$\Phi_1^\eps$, provided that the polarization condition
\eqref{eq:hyppert} is satisfied. By induction, we see that this
perturbs the other $\Phi_j^\eps$'s and ${\bf a}_j^\eps$'s by a
$\O(\delta^\eps)$. Consider the approximate solution defined by 
\begin{equation*}
  u_K^\eps(t,x) = a_0^\eps(x) \exp\(i\sum_{j=1}^K
  t^{2j-1}\Phi_j^\eps(x)/\eps\). 
\end{equation*}
Formally, we have:
\begin{align*}
  {\bf a}^\eps(t,x)e^{i\Phi^\eps(t,x)/\eps} - u_K^\eps(t,x)&=
  \({\bf a}^\eps(t,x)-a_0^\eps(x)\)e^{i\Phi^\eps(t,x)/\eps}\\
  +a_0^\eps(x)& \Big(\exp\Big( i\Phi^\eps(t,x)/\eps\Big) -
  \exp\Big(i\sum_{j=1}^K 
  t^{2j-1}\Phi_j^\eps(x)/\eps\Big)\Big)\\
&= \O(t^2) + \O\Big( \Big(\Phi^\eps(t,x)- \sum_{j=1}^K
  t^{2j-1}\Phi_j^\eps(x)\Big)/\eps\Big)\\
&= \O(t^2) + \O\(t^{2K+1}/\eps\). 
\end{align*}
We infer that the above quantity is small for times such that
$t^\eps\ll 1$ and $t^\eps\ll \eps^{\frac{1}{2K+1}}$. On the other
hand, Proposition~\ref{prop:estprec1} shows that ${\bf a}^\eps
e^{i\Phi^\eps/\eps} $ is a good approximation of $u^\eps$ for times
such that $t^\eps\ll 1$. Therefore, we expect
\begin{equation*}
  \big\|u^\eps(t^\eps,\cdot)-u_K^\eps(t^\eps,\cdot)\big\|_{L^2\cap
  L^\infty}\ll 1 
  \quad \text{ for } t^\eps\ll \eps^{\frac{1}{2K+1}}. 
\end{equation*}
This can be proved by the analysis presented in \S\ref{sec:super}. 
\smallbreak

The case $K=1$ is of special interest. Indeed, the Laplacian plays no
role in the definition of $u_1^\eps$, and we check that it solves:
\begin{equation*}
  i\eps \d_t u_1^\eps = |u_1^\eps|^2 u_1^\eps\quad ;\quad u_{1\mid
  t=0}^\eps= a_0^\eps.
\end{equation*}
This is the solution of the ordinary differential equation considered
in \cite{CCT2} and \cite{BZ}. The above analysis shows that it is a
reasonable approximation of $u^\eps$ for $0<t^\eps\ll \eps^{1/3}$: see
Remark~\ref{rem:laplace}. 
\smallbreak

In view of Theorem~\ref{theo:strong}, we define $u_K^\eps$ and
$v_K^\eps$ in an obvious way, and we have:
\begin{equation*}
  \big\|u^\eps(t^\eps,\cdot)-u_K^\eps(t^\eps,\cdot)\big\|_{L^2\cap
  L^\infty}+ \big\|v^\eps(t^\eps,\cdot)-v_K^\eps(t^\eps,\cdot)\big\|_{L^2\cap
  L^\infty}\ll 1 
  \quad \text{ for } t^\eps\ll \eps^{\frac{1}{2K+1}}. 
\end{equation*}
So to prove Theorem~\ref{theo:strong}, we just have to compare $
u_K^\eps$ and $v_K^\eps$:
\begin{align}
u_K^\eps(t,x) -   v_K^\eps(t,x) = &a_0(x)\exp\Big(i\sum_{j=1}^K
  t^{2j-1}\Phi_j(x)/\eps\Big)\notag\\
&- \widetilde
  a_0^\eps(x)\exp\Big(i\sum_{j=1}^K 
  t^{2j-1}\widetilde\Phi_j^\eps(x)/\eps\Big)\notag\\
=&\(a_0(x)-\widetilde
  a_0^\eps(x)\)\exp\Big(i\sum_{j=1}^K
  t^{2j-1}\Phi_j(x)/\eps\Big)\notag\\
- \widetilde
  a_0^\eps(x)& \Big( \exp\Big(i\sum_{j=1}^K
  t^{2j-1}\Phi_j(x)/\eps\Big)- \exp\Big(i\sum_{j=1}^K 
  t^{2j-1}\widetilde\Phi_j^\eps(x)/\eps\Big)\Big).\label{eq:princ} 
\end{align}
The first term is of order $\delta^\eps$ by assumption. To estimate the second
term, examine:
\begin{equation*}
 \sum_{j=1}^K
  t^{2j-1}\(\Phi_j(x)-\widetilde\Phi_j^\eps(x)\) /\eps .
\end{equation*}
Since we consider times such that $t^\eps\ll 1$ and that we have seen
that $\Phi_j-\widetilde\Phi_j^\eps=\O(\delta^\eps)$ for $j\ge 2$, the
leading order term is simply:
\begin{align*}
 t\(\Phi_1(x)-\widetilde\Phi_1^\eps(x)\) /\eps
 &=\frac{t}{\eps}\(|\widetilde a_0^\eps(x)|^2- |a_0^\eps(x)|^2\)\\
&= \frac{t}{\eps} \( \operatorname{Re}\(a_0-\widetilde
 a_0^\eps)\overline{a_0}\) +\O\((\delta^\eps)^2\)\).
\end{align*}
By assumption, the modulus of \eqref{eq:princ} behaves like
\begin{equation*}
  \left|a_0(x)\sin\(\frac{t\delta^\eps}{\eps}f(x)\)\right|,
\end{equation*}
for some non-trivial function $f$. The conclusion of
Theorem~\ref{theo:strong} follows easily. 

\providecommand{\bysame}{\leavevmode\hbox to3em{\hrulefill}\thinspace}
\providecommand{\MR}{\relax\ifhmode\unskip\space\fi MR }
\providecommand{\MRhref}[2]{%
  \href{http://www.ams.org/mathscinet-getitem?mr=#1}{#2}
}
\providecommand{\href}[2]{#2}

\end{document}